\documentclass[a4paper,11pt,fullpage]{article}
\usepackage[english]{babel}

\usepackage[left=1in,right=1in,top=1in,bottom=1in]{geometry}

\usepackage{amsmath}
\usepackage{amsthm}
\usepackage{amsfonts}
\usepackage{amssymb}
\usepackage{amsxtra}
\usepackage{latexsym}
\usepackage{ifthen}
\usepackage{cite}
\usepackage{esint}
\usepackage[cp1250]{inputenc}

\usepackage{epic}
\usepackage{eepic}

\usepackage[colorlinks=true, urlcolor=blue, citecolor=red, linkcolor=blue]{hyperref}

\DeclareMathOperator*{\osc}{osc}

\numberwithin{equation}{section}
\newtheorem{theorem}{Theorem}[section]
\newtheorem{lemma}{Lemma}[section]
\newtheorem{proposition}{Proposition}[section]
\newtheorem{remark}{Remark}[section]
\newtheorem{definition}{Definition}[section]

\newtheorem{example}{Example}[section]

\def\XXint#1#2#3{{\setbox0=\hbox{$#1{#2#3}{\int}$}
     \vcenter{\hbox{$#2#3$}}\kern-.5\wd0}}

\begin{document}

\title{The weak Harnack inequality for  unbounded minimizers of  elliptic functionals with generalized Orlicz growth}

\author{ Simone Ciani, Eurica Henriques, Igor I. Skrypnik}

  \maketitle

  \begin{abstract}  
  \noindent In this work we prove that the non-negative functions $u \in L^s_{loc}(\Omega)$, for some $s>0$, belonging to the De Giorgi classes 
	\begin{equation}\label{eq0.1}
	\fint\limits_{B_{r(1-\sigma)}(x_{0})} \big|\nabla \big(u-k\big)_{-}\big|^{p}\,  dx \leqslant \frac{c}{\sigma^{q}} \,\Lambda\big(x_{0}, r, k\big)\bigg(\frac{k}{r}\bigg)^{p}\bigg(\frac{\big|B_{r}(x_{0})\cap\big\{u\leqslant k\big\}\big|}{|B_{r}(x_{0})|}\bigg)^{1-\delta},
	\end{equation}
 under proper assumptions on $\Lambda$, satisfy a weak Harnack inequality with a constant  depending  on the $L^s$-norm of $u$. Under suitable assumptions on $\Lambda$, the minimizers of elliptic functionals with generalized Orlicz growth belong to De Giorgi classes satisfying \eqref{eq0.1}; thus this study gives a wider interpretation of Harnack-type estimates derived to double-phase, degenerate double-phase functionals  and functionals with variable exponents.




\end{abstract}

\pagestyle{myheadings} \thispagestyle{plain}
\markboth{Simone Ciani, Eurica Henriques, Igor I. Skrypnik}
{The weak Harnack inequality...}

\section{Introduction and main results}\label{Introduction}

\noindent Relevant developments in the study of regularity of minima of functionals with non-standard growth of $(p, q)$-type were brought to light by the works of Zhikov \cite{Zhi1, Zhi2}, Marcellini \cite{Mar1, Mar2} and Lieberman \cite{Lie}; and a new direction in analysis appears: the theory of Sobolev spaces with a variable exponent, also known as Sobolev-Orlicz spaces. Several other authors joined the venture of better understanding the qualitative properties of functions related to such kind of functionals  (see for instances \cite{Alk, ArrHue, BarColMin1, BarColMin2, BarColMin3, BenHarHasKar, BurSkr, ColMin1, ColMin2, ColMin3, Fan, FanZha, HarHas, HarHasLee, HarHasToi, HarKinLuk, HarKuuLukMarPar, HasOk, LisSkr, MizOhnShi, Ok, RagTac, Sur2, WanLiuZha} and the references there in). The regularity theory for minima of functionals with non-standard growth has been performed over the last three decades for diverse forms of functionals of the type
\begin{equation} \label{Orlicz-functionsl} \int_\Omega \varPhi(x, \nabla u)\ dx, \end{equation} of Orlicz growth, encompassing therefore the double-phase case $\varPhi(x,v)=|v|^p + a(x) |v|^q$, the variable exponent case $\varPhi(x,v)=|v|^{p(x)}$,the double-phase case $\varPhi(x,v)=|v|^{p(x)} +$ $+|v|^{p(x)} a(x)\log(1+ |v|)$, where the exponent function $p(x)$ satisfies a regularity assumption on its modulus of continuity: either the logarithimic condition  
\begin{equation}\label{logcond} \osc\limits_{B_{r}(x_{0})}p(x)\leqslant  \frac{L}{\log\left(\frac{1}{r}\right)},  \quad 0<r<1 \ , \quad 0< L<\infty \ , 
\end{equation}
or a non-logarithimic condition  
\begin{equation}\label{nonlogcond} \osc\limits_{B_{r}(x_{0})}p(x)\leqslant  \frac{\mu(r)}{\log\left(\frac{1}{r}\right)},  \quad \mbox{where} \quad \lim_{r\rightarrow 0} \mu(r)=\infty \ ,  \quad \lim_{r\rightarrow \infty} \frac{\mu(r)}{\log\left(\frac{1}{r}\right)}=0  
\end{equation}
which gives a precise condition for the smoothness of bounded functions in $W^{1,p(x)}(\Omega)$. Zhikov \cite{Zhi3} generalized \eqref{logcond} in the form $ \osc\limits_{B_{r}(x_{0})}p(x)\leqslant  L \frac{\log\left(\log\left(\frac{1}{r}\right)\right)}{\log\left(\frac{1}{r}\right)}$, $ 1<p\leq p(x)$, $0<L<p/n$, 
guaranteeing  the density of smooth functions in the Sobolev space $W^{1,p(x)}(\Omega)$, from which \eqref{nonlogcond} followed.

\vspace{.3cm}

\noindent The quest of deriving Harnack-type estimates met several authors and another interesting amount of approaches and methods. It is well known (and a feature when adopting Moser's or De Giorgi's methods) that the derived estimates are intrinsic in the sense that they depend on the solution or the minimizer - for locally bounded solutions/minimizers, the constant appearing in the Harnack-type estimates depends on this bound (see for instances \cite{Alk} and \cite{BarColMin1} for the variable exponent case under \eqref{logcond} and the double-phase case, respectively). 
So it is somehow expected that, if one only ask the function $u$ to be locally $s$-integrable, {\it i.e.} $u\in L_{loc}^s$, instead of being locally bounded, then the Harnack constants should depend on the $L^s$-norm of $u$. Harnack-type inequalities for unbounded super-solutions of the $p(x)$-Laplacian under the logarithmic condition \eqref{logcond} were proved in \cite{HarKinLuk}, with a constant depending on $||u||_{L^{s}_{loc}(\Omega)}$, $s>0$. This result was generalized in \cite{BenHarHasKar} for unbounded super-solutions to the Euler-Lagrange equation related to the functional \eqref{Orlicz-functionsl}, i.e.
\begin{equation} \label{Orlicz-equation}
    -\mbox{div} \left(\varPhi\left(x,|\nabla u|\right)\dfrac{\nabla u}{|\nabla u|^{2}}\right)=0,
\end{equation} under the ($A 1-s_{\star}$) logarithmic condition
\begin{equation*}
\sup\limits_{x\in B_{r}(x_{0})}\varPhi(x, v) \leqslant C \inf\limits_{x\in B_{r}(x_{0})}\varPhi(x,v) ,\quad r<v\leqslant \frac{1}{|B_{r}(x_{0})|^{\frac{1}{s}}},\quad \text{for any}\quad B_{8r}(x_{0})\subset \Omega .
\end{equation*}
The case of the so-called non-logarithmic conditions differs substantially from the logarithmic one. This can be seen from the fact that, in the logarithmic case, De Giorgi classes  $DG_{\varPhi}(\Omega)$ reduce to the standard and well-studied De Giorgi classes $DG_{p}(\Omega)$ (see \cite{SavSkrYev, ShaSkrYev} and the examples in Section \ref{examples}). Interior continuity, continuity up to the boundary and Harnack's inequality for bounded solutions to the $p(x)$-Laplacian under the non-logarithmic conditions \eqref{nonlogcond} were proved in \cite{AlkKra, AlkSur1, AlkSur2, AlkSur3, Sur1}. These results were generalized in \cite{HadSkrVoi, ShaSkrVoi, SkrVoi, SkrVoi1, SkrYev} for a wider class of elliptic and parabolic equations under the ($\varPhi_{\mu}$) non-logarithmic conditions 
\begin{equation}\label{eq1.3}
\varPhi^{+}_{B_{r}(x_{0})}\left(\frac{v}{r}\right)\leqslant C(K)\,\mu(r) \,\varPhi^{-}_{B_{r}(x_{0})}\left(\frac{v}{r}\right),\quad r<v\leqslant K,\quad \text{for any}\quad B_{8r}(x_{0})\subset \Omega
\end{equation}
with some non-increasing function $\mu(r)\geqslant 1$. The case of unbounded minimizers from the corresponding De Giorgi classes was studied in \cite{SavSkrYev}. Particularly, the weak Harnack inequality was proved under the condition
\begin{equation}\label{eq1.4}
\varPhi^{+}_{B_{r}(x_{0})}\left(\frac{v}{r}\right)\leqslant C \varPhi^{-}_{B_{r}(x_{0})}\left(\frac{v}{r}\right),\quad r<v\leqslant \frac{\lambda(r)}{|B_{r}(x_{0})|^{\frac{1}{s}}},\quad \text{for any}\quad B_{8r}(x_{0})\subset \Omega,
\end{equation}
where $\lambda(r)\in (0,1)$ is a non decreasing function. This study leads to the question of whether a generalization of the De Giorgi classes with Orlicz growth (studied in \cite{SavSkrYev}), can be defined with the aim, on the one hand, to encompass a wider range of functionals and equations; and on the other hand to be defined through a simple energy inequality, that does not refer to some functional $\Phi$, but rather to the scaling between radii and levels to be taken under consideration. This is precisely the scope of our work. 

\vspace{.5cm}
 
\noindent {\textbf{Main Result.}} Our aim is to prove a Weak Harnack inequality for the non-negative minimizers belonging to the De Giorgi classes $DG^{-}_{p, \Lambda}(\Omega),$ where $\Omega$ is a bounded domain in $\mathbb{R}^{n}$ and $\Lambda\big(x_{0}, r, k\big)$ is a positive function satisfying 
\begin{itemize}
\item[($\Lambda_s$)] There exist $s\in (0, \infty]$, a non-decreasing function $\lambda(r) \in (0, 1]$ and a non-increasing function $\mu(r)\geqslant 1$, such that for any $B_{8r}(x_{0})\subset \Omega$, there holds
\begin{equation*}
\Lambda(x_{0}, r , k)\leqslant \mu\left(\left(\frac{r}{k}\right)^{s_\star}\right),\quad \text{for any}\quad  r\leqslant k\leqslant \dfrac{\lambda(r)}{|B_{r}(x_{0})|^{\frac{1}{s}}}, \quad s_\star= \frac{s}{n+s}.
\end{equation*}
\end{itemize}
\noindent In the particular case $s=\infty$, that is when considering bounded minimizers $u$, we set $s_\star=1$.

\noindent For technical reasons, we also assume that the non-decreasing function $\lambda(\cdot)$ satisfies the property
\begin{itemize}
\item[$(\lambda)$] 
for any $0<r \leq \rho$ 
\begin{equation*}
\lambda(\rho)\leqslant \frac{\rho}{r}\,\,\lambda(r)\, .
\end{equation*}
\end{itemize}

\noindent We note that the function $\lambda(r)=(\log\frac{1}{r})^{-\gamma}$, $\gamma>0$ satisfies the above condition, provided that $r$ is small enough.

\begin{definition}\label{def-DG}
Let $s\ge 0$ be a fixed number, and let $\Lambda: \Omega \times \mathbb{R}^+ \times \mathbb{R}^+ \rightarrow \mathbb{R}^+$ be a function satisfying the properties $(\Lambda_s)$-$(\lambda)$ above. We  say that a measurable function $u:\Omega \rightarrow \mathbb{R}$ belongs to the elliptic class $DG^{-}_{p, \Lambda}(\Omega)$ if $u\in W^{1,p}_{loc}(\Omega)$ and there exist numbers  $c$,  $ q\geqslant p>1$, $0\leqslant \delta <\dfrac{p}{n}$  such that for any ball $B_{8r}(x_{0})\subset \Omega$, any  $k>0$ and any $0< \sigma <1$  the following inequalities hold:

\begin{equation}\label{eq1.1}
\fint\limits_{B_{r(1-\sigma)}(x_{0})} \big|\nabla \big(u-k\big)_{-}\big|^{p}\,  dx \leqslant \frac{c}{\sigma^{q}} \,\Lambda\big(x_{0}, r, k\big)\bigg(\frac{k}{r}\bigg)^{p}\bigg(\frac{\big|B_{r}(x_{0})\cap\big\{u\leqslant k\big\}\big|}{|B_{r}(x_{0})|}\bigg)^{1-\delta}, 
\end{equation}
    
\end{definition}

\noindent where $\displaystyle{\fint\limits_{B_{\rho}(x_{0})} g(x)\,dx:=|B_{\rho}(x_{0})|^{-1}\int\limits_{B_{\rho}(x_{0})} g(x)\,dx}$.

\noindent Having set the framework, our main result reads as follows.

\begin{theorem}\label{th1.1}
Let $s\ge 0$ be a fixed number, and let $\Lambda: \Omega \times \mathbb{R}^+ \times \mathbb{R}^+ \rightarrow \mathbb{R}^+$ be a function satisfying the properties $(\Lambda_s)$-$(\lambda)$ above. If $u$ is a non-negative member of $ DG^{-}_{p, \Lambda}(\Omega)\cap L^{s}_{loc}(\Omega)$, then there exist positive real numbers $C_1$, $C_{2}$ and $\beta$, depending only on $n$, $p$, $q$, $c$, $\delta$ and $||u||_{L^{s}_{loc}(\Omega)}$, such that for, all $0<\theta < \dfrac{1}{C_{1}[\mu(\rho)]^{\beta}+\frac{2}{s}}$, there holds
\begin{equation}\label{eq1.2}
\bigg(\fint\limits_{B_{\rho}(x_{0})}\,u^{\theta}\,\,dx\bigg)^{\frac{1}{\theta}}\leqslant \frac{C_2}{\lambda(\rho)}\, e^{C_1[\mu(\rho)]^{\beta}} \Big\{\inf\limits_{B_{\frac{\rho}{2}}(x_{0})}u+\,\,  \rho \Big\},
\end{equation}
provided that $B_{8\rho}(x_{0})\subset \Omega$. 
\end{theorem}

\noindent Our proof relies on De Giorgi's approach, being the main difficulty related to the so-called result on the expansion of positivity. To better understand the choice of our setting, fix a point $x_o \in \Omega$ and a positive number $\rho$ so  that $B_{8\rho}(x_o) \subset \Omega$; and consider that conditions of the type \eqref{eq1.3} are in order. By having information on the measure of the ''positivity set'' of $u$ over the ball $ B_{\rho}(x_{0})$:
\begin{equation*}
\left|B_{\rho}(x_{0})\cap [u\geqslant N]\right|\geqslant \alpha\,\left|B_{\rho}(x_{0})\right| \ , \quad \mbox{for some} \ \rho, N>0 , \quad 0<\alpha<1,
\end{equation*}
we can apply the Local Clustering Lemma, due to DiBenedetto, Gianazza and Vespri \cite{DiBGiaVes} (see also Lemma \ref{subsect2.1} below), to obtain  $\varepsilon_{0} \in (0,1)$, depending only on the known parameters, and a point $y\in B_{\rho}(x_{0})$ such that
\begin{equation*}
\left|B_{r}(y)\cap \left[u\geqslant \frac{N}{2} \right]\right|\geqslant \frac{1}{2}\, \left|B_{r}(y)\right|,\quad r=\varepsilon_{0} \alpha^{2}\dfrac{\rho}{\mu(\rho)}.
\end{equation*}
The presence of the term $\mu(\rho)$ is due to the condition \eqref{eq1.3}. From this, using standard De Giorgi's 
arguments (see for instance \cite{LadUra} Lemma 6.1, or the original \cite{DG}), we inevitably arrive to the estimate
\begin{equation*}
u(x)\geqslant N\, e^{-\gamma\int\limits^{\rho}_{r}[\mu(t)]^{\beta}\frac{dt}{t}} \geqslant N\, e^{\gamma[\mu(r)]^{\beta}\log\frac{\varepsilon_{0}\alpha^{2}}{\mu(\rho)}},\quad x \in B_{\rho}(x_o),
\end{equation*}
with some $\gamma$, $\beta>0$. However the presence of $\alpha$ in the definition of $r$ does not allow us to obtain the weak Harnack inequality, moreover one would need additional and complicated assumptions on $\mu(\cdot)$. In this work, by replacing condition \eqref{eq1.3} by condition ($\Lambda_s$), we are able to prove the weak Harnack inequality \eqref{eq1.2} (see Section \ref{Sec.3} for details).

\begin{remark}
The study of Harnack inequalities for differential equations goes back  to the twentieth century; when considering uniformly elliptic equations a breakthrough on the subject was presented by Moser in \cite{Moser}. By exploring the works of De Giorgi and Moser and by presenting a new Sobolev inequality on minimal surfaces, Bombieri e Giusti \cite{BomGiusti} proved Harnack-type inequalities for elliptic equations on minimal surfaces; while arguing in a different fashion, DiBenedetto and Trudinger \cite{DiBTru} showed that Harnack inequality can be proved directly for functions in the De Giorgi classes. This is precisely the spirit of this work -  to prove that functions belonging to a certain De Giorgi class satisfy a particular weak Harnack estimate and this even when asking for a condition like $(\Lambda_s)$, with a clear dependence on the radius, to be in force.

\noindent We point out two other observations: the first one, is that condition $(\Lambda_s)$ holds for $s=\infty$ therefore this work differs from previous ones and presents a different perspective and approach; the second one, is that for the logarithmic case, {\it i.e.} $\lambda(r)\equiv 1$ and $\mu(r)\equiv 1$, the function $\Lambda\big(x_{0}, r, k\big)$ is constant and therefore, in the bounded case ($s=\infty$) and under logarithmic conditions, Theorem \ref{th1.1} recovers the known result of \cite{DiBTru}.

\end{remark}

\vspace{.5cm}
 
\noindent {\textbf{Structure of the paper.}} In Section \ref{Sec2}, we present a local clustering result due to DiBenedetto, Gianazza and Vespri \cite{DiBGiaVes} and we prove the De Giorgi-type Lemma and a variant result on the expansion of positivity. In Section \ref{Sec.3}, we prove a result on the expansion of positivity that allow us to prove the main result given by Theorem \ref{th1.1}. Finally, in Section \ref{examples} we present several applications of the weak Harnack inequality derived to elliptic functions satisfying generalised Orlicz growth.

\section{Auxiliary Results}\label{Sec2}

In what follows we write $\gamma$ for a positive constant if it can be quantitatively determined a priory in terms of the 
quantities $n$, $p$, $q$, $s$, $c$, $\delta$ and $d:=||u||_{L^{s}_{loc}(\Omega)}$. The generic constant $\gamma$ may change from line to line.

\subsection{ The Local Clustering Lemma }\label{subsect2.1}

\noindent The following lemma is the local clustering lemma due to DiBenedetto, Gianazza and Vespri \cite{DiBGiaVes}.
\begin{lemma}\label{lem2.1}
Let $K_{\rho}(x_{0})$ be a cube in $\mathbb{R}^{n}$ of edge $\rho$ centered at $x_{0}$ and let $u\in W^{1,1}(K_{\rho}(x_{0}))$ satisfies
\begin{equation}\label{eq2.1}
||(u-k)_{-}||_{W^{1,1}(K_{\rho}(x_{0}))} \leqslant \mathcal{K}\,k\,\rho^{n-1}\,\,\,\,\,\,and\,\,\,\,\,\, |K_{\rho}(x_{0}) \cap [u\geqslant k]| \geqslant \alpha |K_{\rho}(x_{0})|,
\end{equation}
for some $\alpha \in (0,1)$, $k\in\mathbb{R}^{+}$ and $\mathcal{K} >0$. Then for any $\eta \in (0,1)$ and any $\nu\in (0,1)$ there exists $y \in K_{\rho}(x_{0})$ and $\varepsilon=\varepsilon(n) \in (0,1)$ such that
\begin{equation}\label{eq2.2}
|K_{r}(y) \cap [u\geqslant \eta\, k]| \geqslant (1-\nu) |K_{r}(y)|,\,\,\, r:=\varepsilon \alpha^{2}\frac{(1-\eta)\nu}{\mathcal{K}}\,\rho.
\end{equation}
\end{lemma}

\subsection{Measurable theoretical results on the expansion of positivity}\label{subsect2.2}

The following results, a De Giorgi-type Lemma and a variant result on the expansion of positivity, are not entirely new but for the sake of completeness (and to allow the text to be self-contained) we decided to present their proofs. In particular, the proof motivates the choice of $\delta < p/N$ in Definition \ref{def-DG}.

\begin{lemma}\label{DeGiorgi}
For $s \ge 0$ and a function $\Lambda: \Omega \times \mathbb{R}^+\times \mathbb{R}^+ \rightarrow \mathbb{R}^+$ satisfying conditions ($\Lambda_s$) and ($\lambda$), we let $u\in DG^{-}_{p,\Lambda}(\Omega)$, $u \geq 0$ and fix $\eta \in [1/2,1)$. Then, there exists $0<\nu_1 <1$, depending on the data and on $\eta$, such that if 
\begin{equation} \label{DG-HP}\big|B_{\rho}(x_o) \cap  [u\leqslant k] \big|\leqslant \nu_{1} \left(\mu\left(\left(\frac{\rho}{2k}\right)^{s_\star}\right) \right)^{-\frac{1}{p/n -\delta}} |B_{\rho}(x_o)|\end{equation}
then, for $0< k \leq \dfrac{\lambda(\rho)}{|B_{\rho}(x_o)|^{1/s}}$ either $(1-\eta) k \leq \rho$ or 
\[ u(x) \geq (1-\eta) k, \qquad x \in B_{\rho/2}(x_o) \ .\]

\end{lemma}

\begin{proof}
Assume $(1-\eta) k> \rho$. Construct the decreasing sequences of radii
$\rho_j=\dfrac{\rho}{2} \left( 1 + \dfrac{1}{2^j} \right)$ and of levels $k_j= (1-\eta) k+\dfrac{\eta k}{2^j}$, for $j=0,1,\cdots$. Observe that, since $u\in DG^{-}_{p,\Lambda}(\Omega)$, $\mu(.)$ is non-increasing and 
\begin{itemize}
\item[(i)] $(1-\eta) k> \rho$ implies $k_j > \rho \geq \rho_j$
\item[(ii)] $k \leq \dfrac{\lambda(\rho)}{|B_{\rho}(x_o)|^{1/s}}$, $\eta \in [1/2,1)$ and $\lambda (\rho) \leq \dfrac{2}{1+ \dfrac{1}{2^j}} \ \lambda(\rho_j)$ leads to 
\[k_j \leq \left(1-\eta +\dfrac{\eta }{2^j} \right)\dfrac{\lambda(\rho)}{|B_{\rho}(x_o)|^{1/s}} \leq \dfrac{\lambda(\rho_j)}{|B_{\rho_j}(x_o)|^{1/s}}\]
\end{itemize}
from \eqref{eq1.1} one first gets
\begin{eqnarray*}
    \int\limits_{B_{\rho_{j+1}}(x_{0})} \big|\nabla \big(u-k_j\big)_{-}\big|^{p}\,  dx  & \leq &  \gamma \, \mu\left(\bigg( \frac{\rho_j}{k_j}\bigg)^{s_\star}\right) \bigg(\frac{k_j}{\rho_j}\bigg)^{p}\bigg(\frac{\big|B_{\rho_j}(x_{0})\cap\big\{u\leqslant k_j\big\}\big|}{|B_{\rho_j}(x_{0})|}\bigg)^{1-\delta} \big|B_{\rho_{j+1}}(x_{0})\big|\\
    & \leq &  \gamma \, \mu\left(\bigg( \frac{\rho}{2k} \bigg)^{s_\star}\right) \bigg(\frac{k}{\rho}\bigg)^{p}\bigg(\frac{\big|B_{\rho_j}(x_{0})\cap\big\{u\leqslant k_j\big\}\big|}{|B_{\rho_j}(x_{0})|}\bigg)^{1-\delta} \big|B_{\rho_{j+1}}(x_{0})\big|
\end{eqnarray*}
and then deduces an algebraic inequality of the type
\[Y_{j+1} \leq \gamma\bigg(\mu\bigg( \frac{\rho}{2k}\bigg)^{s_\star} \bigg) \,  b^j \, Y_j^{1+\alpha}, \qquad \alpha \in (0,1), \ b>1, \quad \]
from which the result follows from a classical iteration argument (see for instance \cite{LadUra}, Lemma 4.7 page 66). This last inequality is obtained by applying  H\"older's inequality together with the Sobolev embedding: let $0 \leq \zeta \leq 1$ be a usual cut-off function between $B_{\rho_{j+1}}(x_0)$ and $B_{\rho_j}(x_0)$, then 
\begin{eqnarray*}
   \left(\frac{\eta k}{2^{j+1}} \right)^p \big|  B_{\rho_{j+1}}(x_{0}) &\cap& [u \leq k_{j+1}]\big|   \leq  \int\limits_{B_{\rho_{j}}(x_{0})}
[\big(u-k_j\big)_{-}{\zeta}]^{p}\,  dx  \\
& \leq &  \gamma \, 
\left(\int\limits_{B_{\rho_{j}}(x_{0})}[\big(u-k_j\big)_{-}\zeta]^{np/(n-p)} \,  dx \right)^{(n-p)/n} \!\! \big|B_{\rho_{j}}(x_{0}) \cap [u \leq k_{j}]\big|^{p/n}\\
& \leq & \gamma \, \int\limits_{B_{\rho_{j}}(x_{0})} \big|\nabla \left(\left(u-k_j\big)_{-}\zeta\right)\right|^{p}\,  dx  \ \big|B_{\rho_{j}}(x_{0}) \cap [u \leq k_{j}]\big|^{p/n}\\
& \leq &\gamma \, \mu\left(\bigg( \frac{\rho}{2k} \bigg)^{s_\star}\right) \bigg(\frac{k}{\rho}\bigg)^{p}\big|B_{\rho_{j}}(x_{0}) \cap [u \leq k_{j}]\big|^{p/n}\\
&  &  \times \left\{ \bigg(\frac{\big|B_{\rho_j}(x_{0})\cap\big\{u\leqslant k_j\big\}\big|}{|B_{\rho_j}(x_{0})|}\bigg)^{1-\delta} \big|B_{\rho_{j+1}}(x_{0})\big| + 2^{pj} \big|B_{\rho_{j}}(x_{0}) \cap [u \leq k_{j}]\big|\right\} 
\end{eqnarray*}
and then 
\[ Y_{j+1}:=\frac{\big|B_{\rho_{j+1}}(x_{0}) \cap [u \leq k_{j+1}\textcolor{blue}{]}\big|}{\big|B_{\rho_{j+1}}(x_{0})\big|} \leq \gamma \, \frac{2^{pj}}{\eta^p}  \, \mu\left(\bigg( \frac{\rho}{2k} \bigg)^{s_\star}\right) \left(\frac{\big|B_{\rho_{j}}(x_{0}) \cap [u \leq k_{j}]\big|}{\big|B_{\rho_{j}}(x_{0})\big|}\right)^{1+p/n-\delta} \ . \]
The proof is complete once we choose $\nu_1= \left(\frac{\gamma \, 2^{\frac{2p}{p/n - \delta}}}{\eta^p}\right)^{-\frac{1}{p/n-\delta}}$.

\end{proof}

\noindent Next result presents the reduction of the measure of the set where $u$ is close to zero.

\begin{lemma}\label{variantexpanposit}
For $s \ge 0$ and a function $\Lambda: \Omega \times \mathbb{R}^+\times \mathbb{R}^+ \rightarrow \mathbb{R}^+$ satisfying conditions ($\Lambda_s$) and ($\lambda$), we let $u\in DG^{-}_{p,\Lambda}(\Omega)$, $u \geq 0$. Assume $0< k \leq \dfrac{\lambda(r)}{|B_{r}(y)|^{1/s}}$ and 
\begin{equation}\label{condalpha}
\big|B_{r}(y) \cap \big\{u\geqslant k\big\}\big|\geqslant \alpha_o |B_{r}(y)| \ , \qquad \mbox{for some} \ \ \alpha_o \in (0,1).
\end{equation}
 Then, for all $\nu \in (0,1)$ and $b \ge 0$, there exist $C_1, C_2, \beta >1$, depending on the data, $\alpha_o$ and $\nu$, such that either 
\begin{equation}\label{condk}
k \leq C_2 \, e^{C_1 \left(\mu\left(\left(\frac{r}{k}\right)^{s_\star}\right)\right)^\beta}\, r  , \quad \beta= \frac{bp+1}{p-1}
\end{equation}
or 
\begin{equation}\label{condset}
 \left| B_r(y)\cap \left[u\leq k \, e^{-C_1 \left(\mu\left(\frac{r}{k}\right)^{s_\star}\right)^\beta}\right] \right|\leq \nu \, \left(\mu\left(\left(\frac{r}{k}\right)^{s_\star}\right)\right)^{-b} \left| B_r(y)\right| . 
\end{equation}

\end{lemma}

\begin{proof}
Assume \eqref{condk} fails to happen, for $C_1, C_2, \beta >1$ to be made explicit. Construct the decreasing sequence 
of levels  $k_j= \dfrac{k}{2^j}$, for $j=0,\cdots, j_\star-1$ (where $j_\star$ is to be chosen) and apply the discrete isoperimetric inequality (see \cite{LadUra} Lemma 3.5 page 55, or the original \cite{DG} Lemma II page 28) to obtain
\begin{equation}\label{isop}
(k_j - k_{j+1}) \big|B_{r}(y) \cap \big[u \leq k_{j+1} \big]\big|  \leq \gamma \, \frac{r^{n+1}}{\big|B_{r}(y) \cap \big[u\geqslant k_j\big]\big|} \, \int\limits_{B_{r}(y) \cap [k_{j+1} <u <k_j]} \big|\nabla u\big| \ .
\end{equation}
For the time being, consider $j_\star$ chosen and $C_1, C_2 $ taken big enough so that
\[ r \leq k_j \leq k \leq \dfrac{\lambda(r)}{|B_{r}(y)|^{1/s}} \ . \]
Recalling \eqref{eq1.1} and by the assumption \eqref{condalpha}, the previous inequality \eqref{isop} together with H\"older's inequality leads to
\[ \big|B_{r}(y) \cap \big[u \leq k_{j+1} \big]\big|^{\frac{p}{p-1}}  \leq \left(\frac{\gamma}{\alpha_o}\right)^{\frac{p}{p-1}} \, \left(\mu\left(\left(\frac{r}{k}\right)^{s_\star}\right)\right)^{\frac{1}{p-1}} \left|B_{r}(y)\right|^{\frac{1}{p-1}} \, \left| B_{r}(y) \cap [k_{j+1} <u <k_j] \right| \ . \]
By summing up this estimate for $j=0, \cdots, j_\star -1$, we arrive at
\[ \big|B_{r}(y) \cap \big[u \leq k_{j_\star} \big]\big| \leq \frac{\gamma}{\alpha_o} \frac{1}{{j_\star}^{\frac{p-1}{p}}} \, \left(\mu\left(\left(\frac{r}{k}\right)^{s_\star}\right) \right)^{\frac{1}{p}} \left|B_{r}(y)\right| \]
and then take
\[j_\star= \left(\frac{\gamma}{\alpha_o \nu} \right)^{\frac{p}{p-1}} \, \left(\mu\left(\left(\frac{r}{k}\right)^{s_\star}\right) \right)^{\frac{bp+1}{p-1}} \ . \]
From this choice of $j_\star$, we choose
\[C_1=\ln 2 \, \left(\frac{\gamma}{\alpha_o \nu} \right)^{\frac{p}{p-1}} ,  \]
and then any constant $C_2 \geq 1$ is good enough so that $k_j \geq r$; in fact
\[ k_j= \frac{k}{2^j} \geq \frac{C_2 \, e^{C_1 \left(\mu\left(\left(\frac{r}{k}\right)^{s_\star}\right)\right)^\beta}\, r }{2^{j_\star}} = C_2 \ r \geq r \ . \]

\end{proof}

\noindent The combination of the two previous lemmas leads to the following result on the expansion of positivity.

\begin{proposition}\label{propexp} For $s \ge 0$ and a function $\Lambda: \Omega \times \mathbb{R}^+\times \mathbb{R}^+ \rightarrow \mathbb{R}^+$ satisfying conditions ($\Lambda_s$) and ($\lambda$), we let $u\in DG^{-}_{p,\Lambda}(\Omega)$, $u \geq 0$. Assume that
\begin{equation}\label{eq2.3}
 0< k \leqslant \frac{\lambda(4r)}{|B_{4r}(y)|^{\frac{1}{s}}},
\end{equation}
and also that, for some $\alpha_{0}\in (0,1)$,
\begin{equation}\label{eq2.4}
\big|B_{r}(y) \cap \big\{u\geqslant k\big\}\big|\geqslant \alpha_{0} |B_{r}(y)|.
\end{equation}
Then there exist $\bar{C_{1}}$, $\bar{C}_{2}>0$, depending only on the data and $\alpha_{0}$, such that either
\begin{equation}\label{eq2.5}
k\,\leqslant \bar{C}_{2}\,e^{\bar{C_{1}}[\mu((r/k)^{s_\star})]^{\bar{\beta}}}\,r,
\end{equation}
or
\begin{equation}\label{eq2.6}
u(x)\geqslant  \frac{k}{2^{j_1+1}} \, e^{-\bar{C_{1}}[\mu((r/k)^{s_\star})]^{\bar{\beta}}},\quad \text{for all}\quad x\in B_{2r}(y),
\end{equation}
where $\bar{\beta}, j_1>1$ are  fixed positive numbers depending only on the data.
\end{proposition}

\begin{proof}
    Assume \eqref{eq2.5} does not hold, for $\bar{C_{1}}$, $\bar{C}_{2}, \bar{\beta} > 0$ to be chosen. Consider the levels 
    \[ k_j=\dfrac{k}{2^j}, \qquad j=j_1, \cdots, j_\star-1 \ , \]
    being $j_1$ the smallest positive integer verifying $j_1 \geq \dfrac{n}{s} \log_2(5/4)$ and $j_\star, \bar{C_{1}}$, $\bar{C}_{2}, \bar{\beta}$ such that $k_j \geq 5r$, and write estimate \eqref{eq1.1} for the truncated functions $(u-k_j)_-$ over the balls $B_{4r}(y) \subset B_{5r}(y)$. This range of $j$ allows us to have

\begin{eqnarray*}
    \int\limits_{B_{4r}(y)} \big|\nabla \big(u-k_j\big)_{-}\big|^{p}\,  dx  & \leq &  \gamma \, \mu\left(\bigg( \frac{5r}{k_j}\bigg)^{s_\star}\right)\bigg(\frac{k_j}{r}\bigg)^{p}\bigg(\frac{\big|B_{5r}(y)\cap\big\{u\leqslant k_j\big\}\big|}{|B_{5r}(y)|}\bigg)^{1-\delta} \big|B_{4r}(y)\big|\\
    & \leq &  \gamma \, \mu\left(\bigg( \frac{r}{k} \bigg)^{s_\star}\right) \bigg(\frac{k_j}{r}\bigg)^{p} \big|B_{4r}(y)\big|
\end{eqnarray*}
    and then, arguing as in the proof of Lemma \ref{variantexpanposit}, one gets 
\begin{equation*} \begin{aligned}
\frac{k_j}{2} \big|B_{4r}(y) &\cap \big[u \leq k_{j+1} \big]\big|   \leq  \gamma \, \frac{r^{n+1}}{\big|B_{4r}(y) \cap \big[u\geqslant k_j\big]\big|} \, \int\limits_{B_{4r}(y) \cap [k_{j+1} <u <k_j]} \big|\nabla u\big|\, dx \\
& \leq  \frac{\gamma}{\alpha_o} \, r \, \left\{\mu \left(\bigg( \frac{r}{k} \bigg)^{s_\star}\right) \,  \bigg(\frac{k_j}{r}\bigg)^{p} \big|B_{4r}(y)\big|\right\}^{\frac{1}{p}} \, 
\big|B_{4r}(y) \cap [k_{j+1} <u <k_j] \big|^{\frac{p-1}{p}}
\end{aligned} \end{equation*}
and finally one arrives at
\[ \big|B_{4r}(y) \cap \big[u \leq k_{j_\star} \big]\big|   \leq \frac{\gamma}{\alpha_o} \, \frac{1}{(j_\star-j_1)^{\frac{p-1}{p}}}  \, \left(\mu\left(\bigg( \frac{r}{k} \bigg)^{s_\star}\right)\right)^{\frac{1}{p}} \, 
\big|B_{4r}(y)  \big| \ . \]
Now we choose
\[\displaystyle{j_\star}= j_1+ \left(\frac{\gamma}{\alpha_o {\nu_1}}\right)^{\frac{p}{p-1}} \, \left[ \mu\left(\bigg( \frac{r}{k} \bigg)^{s_\star}\right)\right]^{\frac{p}{p-1} \left(\frac{1}{p} +\frac{1}{p/n-\delta}\right)}.\]

\noindent Note that by taking $\bar{C_1}=\ln 2 \, \left(\dfrac{\gamma}{\alpha_o \nu_1}\right)^{\frac{p}{p-1}} $, $\bar{C_2} = 2^{j_1 +4}$, $\eta =1/2$ and observing that 
\[\left[\mu\left(\bigg( \frac{r}{k} \bigg)^{s_\star}\right)\right]^{-\frac{1}{p/n-\delta}} \leq 
\left[\mu\left(\bigg( \frac{4r}{2k_{j_\star}} \bigg)^{s_\star}\right)\right]^{-\frac{1}{p/n-\delta}} \] 
we can apply Lemma \ref{DeGiorgi} and conclude
\[u(x) \geq \frac{k}{2^{j_1+1}} \, e^{-\bar{C_{1}}[\mu((r/k)^{s_\star})]^{\bar{\beta}}},\quad \text{for all}\quad x\in B_{2r}(y), \]
being $\bar{\beta}= \dfrac{p}{p-1} \left(\dfrac{1}{p} +\dfrac{1}{p/n-\delta}\right)$. 
\end{proof}


\section{The Weak Harnack Inequality: proof of Theorem \ref{th1.1}}\label{Sec.3}

\subsection{Expansion of Positivity}\label{subsect3.1}

We present and prove the following result from which Theorem \ref{th1.1} easily follows.

\begin{theorem}\label{th3.1}
Let $u\geq 0$, $u \in DG^{-}_{p, \Lambda}(\Omega) \cap  L^{s}_{loc}(\Omega)$, for some $s>0$, and assume conditions ($\Lambda_s$) and ($\lambda$) are fulfilled.  Assume also that, for some $\alpha\in(0,1)$ and $N>0$, there holds
\begin{equation}\label{eq3.1}
\big|B_{\rho}(x_{0})\cap \{u\geqslant N\}\big|\geqslant \alpha |B_{\rho}(x_{0})|, \quad B_{8\rho}(x_{0})\subset \Omega.
\end{equation}
Then there exist $\bar{C}_{3}$, $\bar{C}_{4}$, $\bar{C}_{5}$,  $\tilde{\beta}>1$, depending only on the data,
such that either
\begin{equation}\label{eq3.2}
 \alpha^{\theta_{0}} \lambda(\rho) N\leqslant \bar{C}_{3}\,\,\,e^{\bar{C}_{4}[\mu(\rho)]^{\tilde{\beta}}} \, \rho ,\qquad  \theta_{0} = \frac{1}{s} + 4 + \bar{C_5} \left[\mu(\rho)\right]^{\tilde{\beta}} ,
\end{equation}
or
\begin{equation}\label{eq3.3}
\alpha^{\theta_{0}} \lambda(\rho) N\leqslant  {\bar{C_4}} \, e^{\bar{C}_{4}[\mu(\rho)]^{\tilde{\beta}}} \,\, m_{\frac{\rho}{2}},\qquad m_{\frac{\rho}{2}}=\inf\limits_{B_{\frac{\rho}{2}}(x_{0})} u.
\end{equation}
\end{theorem}

\begin{proof}
We start by assuming that inequality \eqref{eq3.2} is violated, i.e.
\begin{equation}\label{eq3.4}
\alpha^{\theta_{0}} \lambda(\rho) N > \bar{C}_{3}\,\,\,e^{\bar{C}_{4}[\mu(\rho)]^{\tilde{\beta}}} \, \rho.
\end{equation}
Fix $\varepsilon \in(0,1)$, to be chosen, and consider inequality \eqref{eq1.1} written for $k=\varepsilon \lambda(\rho) N$  over the pair of balls $B_{\rho}(x_{0})$ and $B_{2\rho}(x_{0})$:
\begin{equation}\label{eq3.5}
 \fint\limits_{B_{\rho}(x_{0})}|\nabla(u-\varepsilon \lambda(\rho) N)_{-}|^{p}\,dx \leqslant \gamma \bigg(\frac{\varepsilon \lambda(\rho)N}{\rho}\bigg)^{p}\,\,\,\Lambda\big(x_{0}, 2\rho, \varepsilon \lambda(\rho) N\big).
\end{equation}
To estimate the term on the right-hand side of \eqref{eq3.5} we start by noting that
\begin{equation}\label{eq3.6}
 2^{-\frac{n}{s}}\alpha^{\frac{1}{s}}\,N \leqslant \bigg(\fint\limits_{B_{2\rho}(x_{0})} u^{s}\,dx\bigg)^{\frac{1}{s}}\leqslant \frac{d}{|B_{2\rho}(x_{0})|^{\frac{1}{s}}} , 
\end{equation}
and, once we take $\varepsilon=(1+d)^{-1} 2^{-\frac{n}{s}} \alpha^{\frac{1}{s}}$ and use \eqref{eq3.4} and \eqref{eq3.6}, we get
\begin{equation}\label{eq3.7}
2 \rho \leqslant  \varepsilon \lambda(\rho)N
\leqslant \frac{\lambda(\rho)}{|B_{2\rho}(x_{0})|^{\frac{1}{s}}}\leqslant \frac{\lambda(2\rho)}{|B_{2\rho}(x_{0})|^{\frac{1}{s}}} ,
\end{equation}
provided that $\dfrac{\bar{C}_{3}}{2}(1+d)^{-1} 2^{-\frac{n}{s}}\geqslant 1$. Therefore, by condition ($\Lambda_s$) and \eqref{eq3.1}, \eqref{eq3.5}  translates into
\[
\int\limits_{B_{\rho}(x_{0})}| \nabla(u-\varepsilon \lambda(\rho) N)_{-}|\, dx\leqslant \gamma \, \left[\mu\left(\left(\frac{\rho}{\varepsilon \lambda(\rho)N}\right)^{s_\star}\right)\right]^{1/p}\,\varepsilon \lambda(\rho) N\,\rho^{n-1}, \]
and
\[\big|B_{\rho}(x_{0})\cap \{u\geqslant \varepsilon \lambda(\rho) N\}\big|\geqslant \alpha |B_{\rho}(x_{0})|.  \]
Hence, by the local clustering Lemma \ref{lem2.1}, considering $k=\varepsilon\lambda(\rho) N$, $\mathcal{K}=\gamma \, \left[\mu\left(\left(\frac{\rho}{\varepsilon \lambda(\rho)N}\right)^{s_\star}\right)\right]^{1/p}$ and taking $\eta=\nu=\dfrac{1}{2}$, there exists  a point $y\in B_{\rho}(x_{0})$ and $\varepsilon_{0} \in (0,1)$, depending only on the data, such that
\begin{equation}\label{eq3.8}
\big|B_{r}(y)\cap \{u\geqslant \frac{\varepsilon}{2} \lambda(\rho) N \}\big|\geqslant \frac{1}{2}\; |B_{r}(y)|,\quad r=\frac{\varepsilon_{0} \alpha^{2}\rho}{\left(\mu\left(\left(\frac{\rho}{\varepsilon \lambda(\rho)N}\right)^{s_\star}\right)\right)^{1/p}}.
\end{equation}
Set $k_{0}:=\dfrac{\varepsilon r}{2\rho}\lambda(r) N$, then since $\lambda(r)$ is non-decreasing, from \eqref{eq3.8} we obtain
\begin{equation}\label{eq3.9}
\big|B_{r}(y)\cap \{u\geqslant k_{0} \}\big|\geqslant \frac{1}{2} \; |B_{r}(y)|.
\end{equation}
 By  \eqref{eq3.7}, the definition of $r$ and the fact that $\lambda(.)$ is non-decreasing, we have
\begin{equation}\label{eq3.10}
k_{0}\leqslant \varepsilon \lambda(r) N \leqslant \frac{\lambda(r)}{|B_{2\rho}(x_{0})|^{\frac{1}{s}}}\leqslant \frac{\lambda(4r)}{|B_{4r}(y)|^{\frac{1}{s}}}
\end{equation}
and 
\begin{equation*}
 \frac{\varepsilon \lambda(\rho) N}{\rho}\leqslant \rho^{-\frac{1}{s_{\star}}},\quad \text{ hence}\quad \mu(\rho)\geqslant \mu\left(\left(\frac{\rho}{\varepsilon \lambda(\rho)N}\right)^{s_\star}\right) .
\end{equation*}
Our next goal is to apply Proposition \ref{propexp}. For that purpose observe that,  by \eqref{eq3.4}, condition ($\lambda$) and our choice of $k_{0}$, we have
\begin{eqnarray*}
\frac{k_{0}}{r}& = & \frac{\varepsilon}{2\rho} \, \lambda(r) N\geqslant \frac{\varepsilon}{2}\Big(\frac{r}{\rho}\Big)^{2}\frac{\lambda(\rho) N}{r}
\geqslant  \alpha^{-\theta_{0}+\frac{1}{s} +2}\, (1+d)^{-1}\, 2^{-\frac{n}{s} -1} \, \varepsilon_{0} \,  \bar{C_3} \frac{e^{\bar{C_4} \left(\mu\left(\left(\frac{\rho}{\varepsilon \lambda(\rho)N}\right)^{s_\star}\right) \right)^{\tilde{\beta}}}}{\left(\mu\left(\left(\frac{\rho}{\varepsilon \lambda(\rho)N}\right)^{s_\star}\right)\right)^{1/p}} \\
& \geq &  \bar{C_2} \, e^{\bar{C_1} \left(\mu\left(\left(\frac{\rho}{\varepsilon \lambda(\rho)N}\right)^{s_\star}\right) \right)^{\bar{\beta}}} \geq \bar{C_2} \, e^{\bar{C_1} \left(\mu\left(\left(\frac{r}{k_o}\right)^{s_\star}\right) \right)^{\bar{\beta}}}
\end{eqnarray*}
provided that  $\theta_{0} > 2+\dfrac{1}{s}$, $(1+d)^{-1}\, 2^{-\frac{n}{s}-1} \, \varepsilon_{0}\, \bar{C}_{3}\geqslant \bar{C}_{2}$, $\bar{C}_{4}\geqslant \bar{C_{1}}+1$ and $\tilde{\beta}\geqslant \bar{\beta}+ \dfrac{1}{p}$, where  $\bar{C_{1}}$, $\bar{C_{2}}$ and $\bar{\beta} $ are the constants given by Proposition \ref{propexp}. So we can conclude, by \eqref{eq2.6}, 
\begin{equation}\label{eq3.12}
u(x)\geqslant  \frac{k_{0}}{2^{j_1+1}}\, e^{-\bar{C_{1}}\left[\mu\left(\left( r/k_{0}\right)^{s_\star}\right)\right]^{\bar{\beta}}}:=k_{1},\quad \text{for all}\quad x\in B_{2r}(y).
\end{equation}
Repeating this procedure $j$ times, we have
\begin{equation}\label{eq3.13}
u(x)\geqslant \frac{k_{j-1}}{2^{j_1+1}} \, e^{-\bar{C_{1}}\left[\mu\left(\left(2^{j-1}r/k_{j-1}\right)^{s_\star}\right)\right]^{\bar{\beta}}} ,\quad \text{for all}\quad x\in B_{2^{j}r}(y),
\end{equation}
where
\begin{equation*}
k_{i}:=\frac{k_{i-1}}{2^{j_1+1}} \, e^{-\bar{C_{1}}\left[\mu\left(\left(2^{i-1}r/k_{i-1}\right)^{s_\star}\right)\right]^{\bar{\beta}}}< k_{i-1},\quad i=1,...,j,
\end{equation*}
provided that
\begin{equation}\label{eq3.14}
 k_{i}\,\geqslant \bar{C}_{2}\,e^{\bar{C_{1}}\left[\mu\left(\left(2^{i}r/k_{i}\right)^{s_\star}\right)\right]^{\bar{\beta}}} \,2^{i}r ,\quad i=1,...,j-1.
\end{equation}
We now choose $j$ such that $2^{j}r=2\rho$ and then, since $\dfrac{k_{i}}{2^{i}r}\leqslant \dfrac{k_{0}}{r}$, for $i=1,...,j-1$, inequality \eqref{eq3.13} yields
\begin{eqnarray}\label{eq3.15}
u(x) & \geq & \frac{k_{0}}{2^{(j_1+1)j}} \,\,e^{-\bar{C}_{1}\sum\limits_{i=0}^{j-1}\left[\mu\left(\left(2^{i}r/k_{i}\right)^{s_\star}\right)\right]^{\bar{\beta}}} \ , \qquad \forall x \in B_{2\rho}(y) \nonumber\\
& \geqslant & \frac{k_{0}}{2^{(j_1+1)j}} \,\,e^{-\bar{C_{1}}j\left[\mu\left(\left(r/k_{0}\right)^{s_\star}\right)\right]^{\bar{\beta}}}, \qquad  \forall x \in
B_{\frac{\rho}{2}}(x_{0}) \nonumber \\
& \geq & k_{0} \,\,e^{-\tilde{C_{1}}j\left[\mu\left(\left(r/k_{0}\right)^{s_\star}\right)\right]^{\bar{\beta}}}, \quad \tilde{C_1}=\bar{C_1}+ \ln 2(j_1+1)\} \, \qquad  \forall x \in
B_{\frac{\rho}{2}}(x_{0}) \nonumber \\
&\geqslant & \gamma ^{-1}\,e^{-\gamma [\mu(\rho)]^{\tilde{\beta}}}\, \lambda(\rho)\,N\alpha^{\theta_o},\quad x\in B_{\frac{\rho}{2}}(x_{0}) \ . 
\end{eqnarray}

This last inequality was obtained by realizing that 
\begin{eqnarray}\label{estjmu}
j[\mu((r/k_{0})^{s_\star}]^{\bar{\beta}} & = & \left\{1+\frac{1}{p} \log_2 \left(\mu\left(\left(\frac{\rho}{\varepsilon \lambda(\rho)N}\right)^{s_\star}\right)\right)+\log_2\left(\frac{1}{\varepsilon_{0}\alpha^{2}}\right) \right\}[\mu(r/k_{0})]^{\bar{\beta}} \nonumber \\ 
& \leqslant & \left\{1+\log_2\left(\frac{1}{\epsilon_o}\right)+ \frac{1}{p} \right\} \left[\mu(\rho)\right]^{\bar{\beta}+1}-2 \left[\mu(\rho)\right]^{\bar{\beta}}\log_2\left(\alpha\right).
\end{eqnarray}
and taking $\bar{C_5}= 2\tilde{C_1}$ and $\gamma= \max\left\{ \bar{C_3} \epsilon_o^{-1}; \tilde{C_1}\left(1+\log_2\left(\frac{1}{\epsilon_o}\right)+ \frac{1}{p}\right)+\frac{2}{p} \right\} $.

\noindent To complete the proof, it remains to verify \eqref{eq3.14}. Recalling \eqref{eq3.4}, \eqref{eq3.15}, \eqref{estjmu} and our choices of $k_{0}$, $\beta$ and $\theta_{0}$ we have, for any $i=1,...,j$,
\begin{eqnarray}\label{eq3.18}
\frac{k_{i}}{2^{i}r} & = & \frac{k_{0}}{2^{i(j_1+1)} \, 2^{i}r} \, e^{-\bar{C}_{1}\sum\limits_{l=0}^{i-1}[\mu((2^{l}r/k_{l})^{s_\star})]^{\bar{\beta}}}\geqslant \frac{k_{0}}{2\rho}\, e^{-\tilde{C_{1}}j \left[\mu \left(\left(r/k_{0}\right)^{s_\star}\right)\right]^{\bar{\beta}}}
\geqslant\  \frac{\gamma^{-1}}{2} \frac{\lambda(\rho) N \, \alpha^{\theta_{0}}}{\rho}\, e^{-\gamma [\mu(\rho)]^{\tilde{\beta}}}\nonumber \\
& \geqslant & \bar{C}_{3} \, e^{\left(\bar{C}_{4}-2\gamma\right) [\mu(\rho)]^{\tilde{\beta}}} \geqslant 
\bar{C}_{2}\, e^{{\bar{C}_{1}}[\mu\big((2^{i}r/k_{i})\big)]^{\tilde{\beta}}},
\end{eqnarray}
since $\bar{C}_{3} > \bar{C}_{2}$ and once we choose $\bar{C}_{4}\geq 2\gamma + \bar{C}_{1}$, where $\bar{C}_{1}$ and $\bar{C}_{2}$ are the constants obtained in Proposition \ref{propexp}. 
\end{proof}

\subsection{Proof of Theorem \ref{th1.1}}\label{subsect3.2}

The proof of Theorem \ref{th1.1} is almost standard. Let $N$ be a positive real number. For sure 
\[ \frac{\left|B_\rho(x_o) \cap [u\geq N] \right|}{\left|B_\rho(x_o)\right|}  \in (0,1) \]
and therefore, there exists $j \in \mathbb{N}$ such that 
\[\alpha:= \dfrac{1}{2^{j+1}} \leq  \frac{\left|B_\rho(x_o) \cap [u\geq N] \right|}{\left|B_\rho(x_o)\right|} \leq  \dfrac{1}{2^j}\]
By Theorem \ref{th3.1}, there exist positive constants $\bar{C}_{3}$, $\bar{C}_{4}$ and $\tilde{\beta}$, depending only on the data, such that 
\[ \alpha \leq \left(\dfrac{\bar{m}(\rho)}{N}\right)^{\frac{1}{\theta_o}} \ , \qquad \mbox{where} \; \; \bar{m}(\rho):=\dfrac{e^{\bar{C}_{4}[\mu(\rho)]^{\tilde{\beta}}}}{\lambda(\rho)}\,  \,\left( \bar{C_4} \, m_{\frac{\rho}{2}}+ \bar{C_3}\ \rho\right). \]

\noindent Observe that, for $\theta >0$
\begin{eqnarray*}\label{eq3.19}
\fint\limits_{B_{\rho}(x_{0})}\,u^{\theta}\,dx & = &\frac{\theta}{|B_{\rho}(x_{0})|}\int\limits^{\infty}_{0}|B_{\rho}(x_{0})\cap\{u \geq N\}|\,N^{\theta-1}\,dN\leqslant [\bar{m}(\rho)]^{\theta}+ 2 \theta \int\limits^{\infty}_{\bar{m}(\rho)} \alpha \ N^{\theta-1} \, dN\\
& \leq & [\bar{m}(\rho)]^{\theta} + 2\theta[\bar{m}(\rho)]^{\frac{1}{\theta_o}}\int\limits^{\infty}_{\bar{m}(\rho)} \alpha \ N^{\theta-1-\frac{1}{\theta_o}} \, dN =\left(1+\frac{2\theta\theta_{0}}{1-\theta \theta_{0}}\right)[\bar{m}(\rho)]^{\theta}
\end{eqnarray*}
provided that $\theta<\dfrac{1}{\theta_{0}}$. The proof is complete once we choose $\theta\leq \dfrac{1}{\frac{2}{s} + C_1 [\mu(\rho)]^{\beta}} \leq \dfrac{1}{2\theta_{0}}$ 
\begin{eqnarray*}
\left(\fint\limits_{B_{\rho}(x_{0})}\,u^{\theta}\,dx \right)^{\frac{1}{\theta}} & \leq & \left(1+\frac{2\theta\theta_{0}}{1-\theta \theta_{0}}\right)^{\frac{1}{\theta}}\bar{m}(\rho) \leq e^{4\theta_o} \bar{m}(\rho) \\
& \leq & \dfrac{\left(\bar{C_3}+\bar{C_4} \right) e^{4/s+16}}{\lambda(\rho)}\, e^{(\bar{C}_{4}+\bar{C}_{5})[\mu(\rho)]^{\tilde{\beta}}} \,\left( m_{\frac{\rho}{2}}+ \ \rho\right)\\
& \leq & \dfrac{C_2}{\lambda(\rho)}\, e^{C_1 [\mu(\rho)]^{\beta}} \,\left( m_{\frac{\rho}{2}}+ \ \rho\right) \ ,
\end{eqnarray*}
for $C_2:=\left(\bar{C_3}+\bar{C_4} \right) e^{4/s+16} $,  $C_1 = \max \left\{ \bar{C}_{4}+\bar{C}_{5}, 8+2 \bar{C_5} \right\}$ and $\beta= \tilde{\beta}$.

\section{Applications of Theorem \ref{th1.1}}\label{examples}

\noindent In this Section we present detailed examples of integrals of the calculus of variations to which Theorem \ref{th1.1} can be applied. We start by setting the framework, namely by defining the De Giorgi classes in the framework of Sobolev-Orlicz spaces that are associated to the energy inequalities provided by solutions to \eqref{Orlicz-equation}, or minimizers of \eqref{Orlicz-functionsl}.

\vspace{.2cm}

\noindent Consider a nonnegative function $\varPhi: \Omega \times \mathbb{R}^+_0 \rightarrow \mathbb{R}^+_0$ satisfying: 
\begin{itemize}
\item[\hspace{.5cm}$(\varPhi(x,.))$] For all $x \in \Omega$, $v \rightarrow \varPhi(x,v) \ \ \mbox{is increasing} \ , \ \  \displaystyle{\lim_{v\rightarrow 0} \varPhi(x,v)=0  \ , \ \  \lim_{v\rightarrow +\infty}\varPhi(x,v)=+\infty} $; 
\item[($\varPhi_{p,q}$)] there exist $1<p<q$ such that, for all $x\in \Omega$, there holds
\begin{equation*}
\bigg(\frac{u}{v}\bigg)^{p}\leqslant \frac{\varPhi(x,u)}{\varPhi(x,v)}\leqslant \bigg(\frac{u}{v}\bigg)^{q},\quad u\geq v>0.
\end{equation*}
\end{itemize}

\noindent Let $W^{1,\varPhi(\cdot)}(\Omega)$ denote the class of functions $u\in W^{1,1}(\Omega)$ with 
$\int\limits_{\Omega}\varPhi(x, |\nabla u|)\, dx < \infty$. 

\noindent We say that a measurable function $u:\Omega\rightarrow \mathbb{R}^+_0$ belongs to the elliptic class $DG^{-}_{\varPhi}(\Omega)$ if 
\begin{itemize}
\item[(i)]
$u\in W^{1,\varPhi(\cdot)}_{loc}(\Omega)$ 
\item[(ii)] there exist numbers  $c_{1} >0$ and $q>1$ such that\\
\noindent $\forall B_{8r}(x_{0})\subset\Omega, \forall k>0, \forall \sigma\in(0,1)$ 
\begin{equation}\label{eq1.5}
\int\limits_{B_{r(1-\sigma)}(x_{0})}  \varPhi\big(x,|\nabla (u-k)_{-}|\big) \ dx \leqslant  \frac{c_{1}}{\sigma^{q}}\,\int\limits_{B_{r}(x_{0})}  \varPhi\bigg(x,\frac{(u-k)_{-}}{r}\bigg)\  dx.
\end{equation}
\end{itemize}

\noindent In the following example, we cover the case of point-wise non-logarithmic conditions, which goes back to Zhikov's paper \cite{Zhi3}.
\begin{example}\label{ex1.1}
Assume additionally that the function $\varPhi(x,\cdot)$ satisfies
\begin{itemize}
\item[($\varPhi_{\lambda, \mu}$)] For all $r <k\leqslant \dfrac{\lambda(r)}{|B_{r}(x_{0})|^{\frac{1}{s}}}$,
\begin{equation*}
\varPhi^{+}_{B_{r}(x_{0})}\Big(\frac{k}{r}\Big)\leqslant \mu\big((r/k)^{s_{\star}}\big)\,\,\varPhi^{-}_{B_{r}(x_{0})}\Big(\frac{k}{r}\Big),
\quad \text{for any}\quad B_{8r}(x_{0}) \subset \Omega.
\end{equation*}
\end{itemize}
\noindent Here $\varPhi^{+}_{B_{r}(x_{0})}(u):=\sup\limits_{x\in B_{r}(x_{0})} \varPhi(x,u)$, $\varPhi^{-}_{B_{r}(x_{0})}(u):=\inf\limits_{x\in B_{r}(x_{0})} \varPhi(x,u)$, $u>0$.

\noindent Once we consider 
\[\Lambda(x_{0},r, k):=\dfrac{\varPhi^{+}_{B_{r}(x_{0})}\Big(\frac{k}{r}\Big)}{\varPhi^{-}_{B_{r}(x_{0})}\Big(\frac{k}{r}\Big)} \quad \mbox{and} \quad \delta=0 ,\]  
we get  
\[u \in DG_{\varPhi}^{-} (\Omega) \Longrightarrow u \in DG_{p,\Lambda}^{-}(\Omega) \ .\]
Although the proof is neither new nor difficult, we decided to presented it just to keep the text self-contained. First we need the analogue of Young's inequality: being $\varphi_{\gamma}(x,v):=\dfrac{\varPhi(x,v)}{v^{\gamma}}$ and $v>0$,
\begin{equation*}
\varphi_{p}(x,a) b^{p} \leqslant \varPhi(x,a) + \varPhi(x,b),\quad a,b >0.
\end{equation*}
Indeed, if $b\leqslant a$ then $\varphi_{p}(x,a) b^{p} \leqslant \varphi_{p}(x,a) a^{p}=\varPhi(x,a)$ and if $b\geqslant a$, using the fact that  the function $\varphi_{p}(x,a)$ is increasing, we obtain  $\varphi_{p}(x,a) b^{p} \leqslant \varphi_{p}(x,b) b^{p}=\varPhi(x,b)$.

\noindent Assume that $u \in DG_{\varPhi}^{-} (\Omega)$. By \eqref{eq1.5} we have
\[\hspace{-4cm}\varPhi^{-}_{B_{r}(x_{0})}\left(\frac{k}{r}\right)\,\, \left(\frac{r}{k}\right)^p \, \int\limits_{B_{(1-\sigma)r}(x_{0})} |\nabla (u-k)_{-}|^{p}\,\,dx \]
\begin{eqnarray*}
& \leqslant & 
\int\limits_{B_{(1-\sigma)r}(x_{0})}\varphi_{p}\left(x, \frac{k}{r}\right) |\nabla (u-k)_{-}|^{p}\,\,dx \\
& \leqslant &  \int\limits_{B_{(1-\sigma)r}(x_{0})}\varPhi\left(x,|\nabla(u-k)_{-}|\right)\,\,dx+\int\limits_{B_{(1-\sigma)r}(x_{0})\cap [u\leqslant k]}\varPhi\left(x,\frac{k}{r}\right)\,\,dx \\
& \leqslant &
\frac{c_{1}+1}{\sigma^{q}}\,  \varPhi^{+}_{B_{r}(x_{0})}\left(\frac{k}{r}\right) \left|B_{(1-\sigma)r}(x_{0})\cap\big[u\leqslant k\big]\right|,
\end{eqnarray*}
which yields
\begin{equation*}
\int\limits_{B_{(1-\sigma)r}(x_{0})} |\nabla (u-k)_{-}|^{p}\,\,dx \leqslant \frac{c_{1}+1}{\sigma^{q}}\,\frac{\varPhi^{+}_{B_{r}(x_{0})}\Big(\frac{k}{r}\Big)}{\varPhi^{-}_{B_{r}(x_{0})}\Big(\frac{k}{r}\Big)} \Big(\frac{k}{r}\Big)^{p}\,\big|B_{r}(x_{0})\cap\big[u\leqslant k\big]\big|
\end{equation*}
and therefore $u \in DG_{p,\Lambda}^{-}(\Omega)$.

\vspace{.3cm} 

\noindent In what follows we present three specific functions $\varPhi_i(x, .)$, $i=1,2,3$, that verify assumption ($\varPhi_{\lambda, \mu}$), hence one can say that any minimizer $u$ for the functionals related to $\varPhi_i(x, .)$, $i=1,2,3$, belongs to $DG_{p,\Lambda}^{-}(\Omega)$ (and Theorem \ref{th1.1} applies).

\begin{itemize}
    \item[(1)]  The function $\varPhi_{1}(x, v)=v^{p}+a(x) v^{q}$, $a\geq 0$, satisfies condition ($\varPhi_{\lambda, \mu}$) with 
    \[ \lambda(r)=[\bar{\mu}(r)]^{-\frac{1}{q-p}}, \quad \mu(\cdot)\equiv const \quad \mbox{and} \quad s= \dfrac{(q-p)n}{\alpha+p-q}, \] provided that 
    \begin{equation} \label{conda} \osc\limits_{B_{r}(x_{0})}a(x)\leqslant A r^{\alpha} \bar{\mu}(r) \ , \quad  A>0 \quad \mbox{and} \quad \alpha >q-p .
    \end{equation}
In fact, let $\varPhi^{+}_1\Big(\frac{k}{r}\Big):= \sup\limits_{x\in B_{r}(x_{0})} \varPhi_1\Big(x, \frac{k}{r}\Big)$ and $\varPhi^{-}_{1}\Big(\frac{k}{r}\Big):=\inf\limits_{x\in B_{r}(x_{0})} \varPhi_1\Big(x, \frac{k}{r}\Big)$.
    \begin{eqnarray*} 
   \varPhi^{+}_1\Big(\frac{k}{r}\Big) & \leqslant & \varPhi^{-}_{1}\Big(\frac{k}{r}\Big) + A r^\alpha \bar{\mu}(r) \Big(\frac{k}{r}\Big)^q \\
    & \leq & \left( 1 + Ar^\alpha \bar{\mu}(r) \Big(\frac{k}{r}\Big)^{q-p} \right)\varPhi^{-}_{1}\Big(\frac{k}{r}\Big)\\ 
    & \leq & \left( 1 + A \chi_n^{(p-q)/s}\, r^{\alpha -q+p -\frac{n}{s}(q-p)}\bar{\mu}(r) \lambda(r)^{q-p}\right)\varPhi^{-}_{1}\Big(\frac{k}{r}\Big) \\
    & \leq  & C(A,p,q,n, s,\alpha,\Omega) \, \varPhi^{-}_{1}\Big(\frac{k}{r}\Big)\ .
    \end{eqnarray*}

    \item[(2)] The function $\varPhi_{2}(x,v)=v^{p}+a(x) v^{q}\Big(1+\log(1+v)\Big)^{L_{1}}$, with $L_{1}\geqslant 0$, satisfies condition ($\varPhi_{\lambda, \mu}$) with 
    \[\lambda(r)=[\bar{\mu}(r)]^{-\frac{1}{q-p}}\left[\log\left(\dfrac{1}{r}\right) \right]^{-\frac{L_{1}}{q-p}}, \quad \mu(\cdot)\equiv const \quad \mbox{ and} \quad s\geqslant \dfrac{(q-p)n}{\alpha+p-q} \] provided that \eqref{conda} holds (the proof follows closely the arguments presented previously).

    \item[(3)] The function $\varPhi_{3}(x,v)=v^{p(x)}$, $\osc\limits_{B_{r}(x_{0})}p(x)\leqslant L\, \dfrac{\log\left(\bar{\mu}(r)\right)}{\log\left(\frac{1}{r}\right)}$, $L>0$, satisfies condition ($\varPhi_{\lambda, \mu}$) with $\lambda(\cdot)\equiv 1$, 
$\mu(\cdot)=[\bar{\mu}(\cdot)]^{Ls_\star}$, provided that  $f(r):= \dfrac{\log\left(\bar{\mu}(r)\right)}{\log\left(\frac{1}{r}\right)}$ is non-decreasing.

Let $\varPhi^{+}_3\Big(\frac{k}{r}\Big):= \sup\limits_{x\in B_{r}(x_{0})} \varPhi_3\Big(x, \frac{k}{r}\Big)$, $\varPhi^{-}_{3}\Big(\frac{k}{r}\Big):=\inf\limits_{x\in B_{r}(x_{0})} \varPhi_3\Big(x, \frac{k}{r}\Big)$ and take $\lambda(r)=1$. 

\noindent Observe that, since
\begin{eqnarray*}
 r < k\leqslant \dfrac{1}{|B_{r}(x_{0})|^{\frac{1}{s}}} & \Longrightarrow & \left\{ \begin{array}{l}
 \dfrac{k}{r} >1 \\
 r \leq \left(\dfrac{r}{k}\right)^{s_\star}
 \Longrightarrow 
 f(r) \leq f\left( \left(\dfrac{r}{k}\right)^{s_\star}\right) 
 \end{array}\right.
\end{eqnarray*}
we get
\[ \log \left( \left(\frac{k}{r}\right)^{\dfrac{\log\left(\bar{\mu}(r)\right)}{\log\left(\frac{1}{r}\right)}}\right) = \log\left(\bar{\mu}(r)\right) \ \frac{\log\left(\frac{k}{r}\right)}{\log\left(\frac{1}{r}\right)} \leq \log \left( \bar{\mu}((r/k)^{s_\star})\right)^{s_\star}
\]
and then 
\begin{eqnarray*} 
      \varPhi^{+}_3\left(\frac{k}{r}\right) & = & \left(\frac{k}{r}\right)^{\osc\limits_{B_{r}(x_{0})}p(x)} \,  \varPhi^{-}_3\Big(\frac{k}
     {r}\Big)\leq \left(\frac{k}{r}\right)^{L\dfrac{\bar{\mu}(r)}{\log\left(\frac{1}{r}\right)}} \,  \varPhi^{-}_3\Big(\frac{k}
     {r}\Big) \\
     & \leq & \left( \bar{\mu}((r/k)^{s_\star})\right)^{Ls_\star}\,  \varPhi^{-}_3\Big(\frac{k}
     {r}\Big) 
    \end{eqnarray*}

\end{itemize}

\noindent Accordingly, for example, condition ($\varPhi_{\lambda, \mu}$) holds for $\varPhi_{1}(x, \cdot)$ and $\varPhi_{2}(x, \cdot)$ with $p=p(x)$ and/or $q=q(x)$ for an appropriate choice of the functions $\lambda(\cdot)$, $\mu(\cdot)$.

\begin{remark}
In \cite{SavSkrYev} one can find Harnack inequalities for minimizers of non-uniformly elliptic functionals that belong to a De Giorgi class $u \in DG_{\varPhi}^{-}(\Omega)$ under a setting quite close to this one. 
\end{remark}

\end{example}

\noindent In the following example we  cover the non-uniformly elliptic case, which goes back to Trudinger's paper \cite{Tru} (see also \cite{BelSch}).
\begin{example}\label{ex1.2}
Consider that the function $\varPhi(x,\cdot)$ satisfies the extra condition
\begin{itemize}
\item[($\varPhi^{\prime}_{\lambda, \mu}$)] For all $ r <k\leqslant \dfrac{\lambda(r)}{|B_{r}(x_{0})|^{\frac{1}{s}}}$ there holds
\begin{equation*}
\Lambda_{+,\varPhi}\Big(x_{0}, r, \frac{k}{r}\Big) \Lambda_{-,\varPhi}\Big(x_{0}, r, \frac{k}{r}\Big)\leqslant \mu(r),
\quad \text{for any}\quad B_{8r}(x_{0}) \subset \Omega,
\end{equation*}
\end{itemize}
\noindent being $u>0$, $m,p>1$, $t>\max\{1, p-1\}$ and such that $\dfrac{1}{m}+\dfrac{1}{t}<\dfrac{p}{n}$, and 
\[ \Lambda_{+,\varPhi}\big(x_{0}, r, u\big):=\Big(\fint\limits_{B_{r}(x_{0})}[\varPhi(x,u)]^{m}\,dx\Big)^{\frac{1}{m}}, \quad  \Lambda_{-,\varPhi}\big(x_{0}, r, u\big):=\Big(\fint\limits_{B_{r}(x_{0})}[\varPhi(x,u)]^{-t}\,dx\Big)^{\frac{1}{t}} \ . \]

\noindent For 
\[ \Lambda(x_{0},r,k)=\left[\Lambda_{+,\varPhi}\left(x_{0}, r, \frac{k}{r}\right) \Lambda_{-,\varPhi}\left(x_{0}, r, \frac{k}{r}\right)\right]^{\frac{t}{t+1}} \quad \mbox{and} \quad \delta= \dfrac{t}{t+1}\left(\dfrac{1}{t}+\dfrac{1}{m}\right) , \]  
we get  
\[u \in DG_{\varPhi}^{-} (\Omega) \Longrightarrow u \in DG_{\frac{pt}{t+1},\Lambda}^{-}(\Omega) \ .\]

\noindent We show this simple fact for  completeness of the proof. By H\"{o}lder's and  Young's inequalities
\[
\fint\limits_{B_{r}(x_{0})}|\nabla(u-k)_{-}|^{\frac{pt}{t+1}}\,dx  \leqslant \left(\fint\limits_{B_{r}(x_{0})}\left[\varphi_{p}\left(x,\frac{k}{r}\right)\right]^{-t}dx\right)^{\frac{1}{t+1}}\left(\fint\limits_{B_{r}(x_{0})}\varphi_{p}\left(x,\frac{k}{r}\right)|\nabla(u-k)_{-}|^{p}dx\right)^{\frac{t}{t+1}} \]
\begin{eqnarray*}
& \leq &  \left(\fint\limits_{B_{r}(x_{0})}\left[\varphi_{p}\left(x,\frac{k}{r}\right)\right]^{-t}\,dx\right)^{\frac{1}{t+1}} \left(\fint\limits_{B_{r}(x_{0})}\varPhi\left(x,\frac{k}{r}\right)\,dx+\fint\limits_{B_{r}(x_{0})}\varPhi(x,|\nabla (u-k)_{-}|)\,dx\right)^{\frac{t}{t+1}} \\ 
& \leqslant &  \frac{c}{\sigma^{\frac{qt}{t+1}}}
\, \left(\frac{k}{r}\right)^{\frac{pt}{t+1}}\, \left(\fint\limits_{B_{r}(x_{0})}\left[\varPhi\left(x,\frac{k}{r}\right)\right]^{-t}\,dx\right)^{\frac{1}{t+1}}\left(\fint\limits_{B_{r}(x_{0})}\left[\varPhi\left(x,\frac{k}{r}\right)\right]^{m}\,dx\right)^{\frac{t}{m(t+1)}}\\
& & \times 
\bigg(\frac{|B_{r}(x_{0})\cap\{u\leqslant k\}|}{|B_{r}(x_{0})|}\bigg)^{\frac{(m-1)t}{m(t+1)}},
\end{eqnarray*}
which yields
\begin{equation*}
\fint\limits_{B_{(1-\sigma)r}(x_{0})}|\nabla(u-k)_{-}|^{\frac{pt}{t+1}}\,dx\leqslant\frac{c}{\sigma^{\frac{qt}{t+1}}}\, 
\Lambda\big(x_{0},r,k\big) \, \left(\frac{k}{r}\right)^{\frac{pt}{t+1}}\left(\frac{|B_{r}(x_{0})\cap\{u\leqslant k\}|}{|B_{r}(x_{0})|}\right)^{1-\frac{t}{t+1}(\frac{1}{t}+\frac{1}{m})}.
\end{equation*}

\noindent As an example, consider the function $\varPhi_{5}(x, v)=v^{p(x)}$, $p(x)=p+L\dfrac{\log\left(\log\left(\frac{1}{|x-x_o|}\right)\right)}{\log\left(\frac{1}{|x-x_o|}\right)}$. This function 
 satisfies condition ($\varPhi^{'}_{\lambda,\mu}$) with $\lambda(\cdot)=\mu(\cdot)\equiv 1$, provided that $r$ is small enough (see
\cite{SkrYev} for details).

\end{example}

\vskip0.4cm
\noindent{\bf Acknowledgements.} S. Ciani acknowledges the support of the department of Mathematics of the University of Bologna
Alma Mater, and of the PNR fundings 2021-2027. E. Henriques was financed by Portuguese Funds through
FCT - Funda\c c\~ao para a Ci\^encia e a Tecnologia - within the Projects UIDB/00013/2020 and
UIDP/00013/2020. I. Skrypnik is partial supported by a grant from the Simons Foundation (Award 1160640, Presidential Discretionary-Ukraine Support Grants, Skrypnik I.I.).

\vskip3.5mm

\noindent{\bf Research Data Policy and Data Availability Statements.}  All data generated or analysed
during this study are included in this article.

\bigskip

CONTACT INFORMATION

\noindent \textbf{Simone Ciani}\\
Department of Mathematics of the University of Bologna, Piazza Porta San Donato, 5, 40126 Bologna, Italy \\
 simone.ciani3@unibo.it

\medskip

\noindent \textbf{Eurica Henriques}\\
Centro de Matem\'atica, Universidade do Minho - Polo CMAT-UTAD \\
Departamento de Matem\'atica - Universidade de Tr\'as-os-Montes e Alto Douro, 5000-801 Vila Real, Portugal\\
eurica@utad.pt

\medskip
\noindent \textbf{Igor I.~Skrypnik}\\Institute of Applied Mathematics and Mechanics,
National Academy of Sciences of Ukraine, \\ 
 Batiouk Str. 19, 84116 Sloviansk, Ukraine\\
Vasyl' Stus Donetsk National University, 600-richcha Str. 21, 21021 Vinnytsia, Ukraine\\ihor.skrypnik@gmail.com
\end{document}